\documentclass[10pt,leqno]{article}
\usepackage{amssymb,amsbsy,amsmath,amsfonts,amssymb,amscd, mathrsfs}

\usepackage[english]{babel}
\usepackage[T1]{fontenc}
\usepackage{indentfirst}

\makeatletter
\@addtoreset{equation}{section}
\makeatother

\newtheorem{statement}{}[section]
\newtheorem{theoreme}[statement]{Theorem}
\newtheorem{lemme}[statement]{Lemma}

\newtheorem{proposition}[statement]{Proposition}

\newtheorem{corollaire}[statement]{Corollary}

\newcommand\C{\mathbb C}

\newcommand\T{\mathbb T}
\newcommand\D{\mathbb D}
\newcommand\Z{\mathbb Z}
\newcommand\e{{\rm e}}

\newcommand\eps{\varepsilon}
\newcommand\ind{{\rm 1\kern-.30em I}}
\newcommand\qed{\hfill $\square$}
\renewcommand \Re{{\mathfrak R}{\rm e}\,}
\renewcommand \Im{{\mathfrak I}{\rm m}\,}

\title{\bf Some examples of compact composition operators on $H^2$}
\author{\it Pascal Lef\`evre, Daniel Li,\\ \it Herv\'e Queff\'elec, Luis Rodr{\'\i}guez-Piazza}

\date{\footnotesize \today}

\begin{document}

\maketitle

\noindent{\bf Abstract.} \emph{We construct, in an essentially explicit way, various composition 
operators on $H^2$ and study their compactness or their membership in the Schatten classes. We construct:  
non-compact composition operators on $H^2$ whose symbols have the same modulus on the boundary of 
$\D$ as symbols whose composition operators are in various Schatten classes $S_p$ with $p > 2$; compact 
composition operators which are in no Schatten class but whose symbols have the same modulus on the 
boundary of $\D$ as symbols whose associated composition operators are in $S_p$ for every $p>2$.}
\medskip

\medskip


%%%%%%%%%%%%%%

\section{Introduction}

Compactness of composition operators on $H^2$ was first studied in 1968 by H. Schwartz in his 
doctoral dissertation \cite{Schw}, and refined in 1973 by J.~Shapiro and P. Taylor \cite{Shap-Tay}, who 
discovered the role played by the classical angular derivative, and refined the compactness problem by 
asking which composition operators belong to various Schatten classes $S_p$. In particular, they showed 
in \cite{Shap-Tay} that $C_\phi \in S_2$, the Hilbert-Schmidt class, if and only if 
$\int_{\partial \D} (1 - |\phi^\ast|)^{-1} < +\infty$, where $\phi^\ast$ denotes the radial limit function 
of $\phi$. In this paper, we show that for the larger class $S_p$ with $p > 2$, the situation is 
completely different: we prove in Theorem~\ref{meme module Schatten}, that for every $p > 2$, there exist 
two symbols $\phi_1$ and $\phi_2$ having the same modulus on $\partial \D$ and such that $C_{\phi_1}$ 
is not compact on $H^2$, but $C_{\phi_2}$ is in the Schatten class $S_p$.\par
\medskip

An amusing feature of the theory of composition operators is that, whereas sophisticated necessary and 
sufficient conditions for the composition operators $C_\phi \colon H^2 \to H^2$ to belong to the  
Schatten classes $S_p = S_p (H^2)$ have been known for more than twenty years (\cite{Lue}, \cite{Lue-Zhu}), 
either in terms of the Nevanlinna counting function, or in terms of the pull-back measure $m_\phi$,  yet 
explicit and concrete examples are lacking. For example, D. Sarason posed in 1988 the question 
whether there existed a compact composition operator $C_\phi \colon H^2 \to H^2$ which was in \emph{no} 
Schatten class, and the (affirmative) answer only came in 1991, by T. Carroll and C. Cowen (\cite{Carroll-Cowen}). 
Their example, based on the Riemann mapping Theorem was not completely explicit. Moreover, their 
construction used a difficult and delicate argument, with estimates of the hyperbolic metric for certain domains,  
due to Hayman (see however \cite{Zhu} and \cite{Jones}). We shall see in this paper that Luecking's criterion
\cite{Lue} for pullback measures leads to very concrete examples of composition operators in various Schatten 
classes.\par

In Section~\ref{Preliminaries}, we give a necessary condition, 
Proposition~\ref{condition necessaire Schatten}, on the Carleson function $\rho_\phi$ in order for the 
composition operator $C_\phi$ to be in $S_p$, as well as a general construction of symbols $ \phi$ with 
control on their Carleson function.\par

In Section~\ref{section same modulus}, we construct, for every $p > 2$, symbols $\phi_1$ and $\phi_2$ 
having the same modulus on $\partial \D$ such that $C_{\phi_1}$ is not compact on $H^2$, but 
$C_{\phi_2}$ is in the Schatten class $S_p$ (Theorem~\ref{meme module Schatten}). \par

In Section~\ref{log puissance theta}, we revisit an example of J. Shapiro and P. Taylor (\cite{Shap-Tay}, \S~4) 
to show that for every $p_0 > 0$, there exists a symbol  $\phi $ such that the composition 
operator $C_\phi \colon H^2 \to H^2$ is in the Schatten class $S_p$ for every $p > p_0$, but not in 
$S_{p_0}$ (Theorem~\ref{theo Schatten}), and also that for every $p_0 >0$, there exists a symbol $\phi$ 
such that $C_\phi \colon H^2 \to H^2$ is in the Schatten class $S_{p_0}$, but not in $S_p$, for $p < p_0$ 
(Theorem~\ref{theo Schatten bis}). Moreover, there exists a symbol $\phi$ such that 
$C_\phi \colon H^2 \to H^2$ is compact but in no Schatten class $S_p$ for $p < \infty$ 
(Theorem~\ref{no Schatten}) and there exists a symbol $\psi$, whose boundary values $\psi^\ast$ have the 
same modulus as those $\phi^\ast$ of $\phi$ on $\partial \D$, but for which $C_\psi \colon H^2 \to H^2$ is 
in $S_p$ for every $p>2$ (Theorem~\ref{no Schatten meme module}).\par 
After our work was completed, we became aware of the papers \cite{Zhu} and \cite{Jones}; in \cite{Jones}, 
M. Jones gives another proof of the theorem of Carroll and Cowen (our Theorem~\ref{no Schatten}), and Y. Zhu 
gives also another proof of this theorem, as well as a proof of our Theorem~\ref{theo Schatten}. However, our 
proofs are different and lead to further results: see Theorem~\ref{meme module Schatten}.\par 
\bigskip

\noindent{\bf Acknowledgements.} We thanks the referee for a very careful reading of this paper and 
many suggestions to improve its writing.

\section{Notation}

Throughout this paper, the notation $f\approx g$ will mean that there are two constants $0 < c < C < +\infty$ 
such that $ c f(t) \leq g(t) \leq C f(t)$ (for $t$ sufficiently near of a specified value), and the notation 
$ f (t) \lesssim g (t)$, when $t$ is in the neighbourhood of some value $t_0$, will have the same meaning as  
$g = O (f)$.\par
\medskip

We shall denote by $\D$ the open unit disc of the complex plane: $\D = \{ z\in \C\,; \ |z| < 1\}$, and by 
$\T = \partial \D$ its boundary: $\T =\{ z\in \C\,;\ |z | = 1\}$. We shall denote by $m$ the normalized 
Lebesgue measure on $\T$.\par
For every analytic self-map $\phi \colon \D \to \D$, the composition operator $C_\phi$ is the map 
$f \mapsto f\circ \phi$. By Littlewood's subordination principle (see \cite{Duren}, Theorem~1.7), every 
composition operator maps every Hardy space $H^p$ ($p >0$) into itself, and is continuous on $H^p$.\par
\smallskip

For every $\xi \in \T$ and $0 < h < 1$, the Carleson window $W (\xi, h)$ is the set
\begin{displaymath}
W (\xi, h) = \{ z \in \D\,;\ |z| \geq 1 - h\quad \text{and} \quad |\arg ( z \bar{\xi})| \leq h\}.
\end{displaymath}
For every finite positive measure $\mu$ on $\D$, one sets:
\begin{displaymath}
\rho_\mu (h) = \sup_{\xi \in \T} \mu[W (\xi, h)].
\end{displaymath}
We shall call this function $\rho_\mu$ the \emph{Carleson function} of $\mu$.\par
When $\phi \colon \D \to \D$ is an analytic self-map of $\D$, and $\mu =m_\phi$ is the 
measure defined on $\D$, for every Borel set $B \subseteq \D$, by:
\begin{displaymath}
m_\phi (B) = m (\{\xi\in \T\,;\ \phi^\ast (\xi) \in B \}),
\end{displaymath}
where $\phi^\ast$ is the boundary values function of $\phi$, we shall denote $\rho_{m_\phi}$ by 
$\rho_\phi$. In this case, we shall say that $\rho_\phi$ is the Carleson function of $\phi$.\par
\smallskip

For $\alpha \geq 1$, we shall say that $\mu$ is an \emph{$\alpha$-Carleson measure} if 
$\rho_\mu (h) \lesssim h^{\alpha}$. For $\alpha =1$, $\mu$ is merely said to be a \emph{Carleson measure}. 
\par\smallskip

The Carleson Theorem (see \cite{Duren}, Theorem~9.3) asserts that, for $1 \leq p < \infty$ (actually, for 
$ 0 < p < \infty$), the canonical inclusion $j_\mu \colon H^p \to L^p (\mu)$ is bounded if and only if 
$\mu$ is a Carleson measure. Since every composition operator $C_\phi$ is continuous on $H^p$, it defines a 
continuous map $j_\phi \colon H^p \to L^p (\mu_\phi)$; hence every pull-back measure $\mu_\phi$ is a 
Carleson measure.\par
When $C_\phi \colon H^2 \to H^2$ is compact, one has, as it is easy to see:
\begin{equation}\label{condition phi plus petit que un au bord}
|\phi^\ast| < 1 \qquad \text{\emph{a.e.} on } \partial \D.
\end{equation}
Hence, we shall only consider in this paper symbols $\phi$ for which 
\eqref{condition phi plus petit que un au bord} is satisfied (which is the case, as we said, when $C_\phi$ is 
compact on $H^2$).\par
\smallskip

B. MacCluer (\cite{McCluer}, Theorem~1.1) has shown (assuming 
condition~\eqref{condition phi plus petit que un au bord}) that $C_\phi$ is compact on $H^p$ if and 
only if $\rho_\phi (h) = o\, (h)$, as $h$ goes to $0$.
\medskip

Note that, in this paper, we shall not work, most often, with exact inequalities, but with 
inequalities up to constants. It follows that we shall not actually work with true Carleson windows 
$W (\xi, h)$ (or Luecking sets, defined below), but with distorted Carleson windows: 
\begin{displaymath}
\tilde W (\xi, h) =\{ z\in \D\,;\ |z| \geq 1 - ah \quad \text{and} \quad | \arg (z \bar{\xi} )|\leq bh\},
\end{displaymath}
where $a,b > 0$ are given constants. Since, for a given symbol $\phi$, one has:
\begin{displaymath}
m_\phi \big(W (\xi, c\,h) \big) \leq m_\phi \big(\tilde W (\xi, h) \big) \leq m_\phi \big(W (\xi, C\,h) \big) 
\end{displaymath}
for some constants $c= c(a,b)$ and $C = C (a,b)$ which only depend on $a$ and $b$, that will not matter 
for our purpose.

\section{Preliminaries}\label{Preliminaries}

\subsection{Luecking sets and Carleson windows}

We shall begin by recalling the characterization, due to D. Luecking (\cite{Lue}), of the composition operators 
on $H^2$ which belong to the Schatten classes. Let, for every integer $n \geq 1$ and $0 \leq j \leq 2^n -1$: 
\begin{displaymath}
R_{n,j} = \Big\{z\in \D\,;\ 1 - 2^{-n} \leq |z| < 1 - 2^{-n - 1}\quad \text{and}\quad 
\frac{2j\pi}{2^n} \leq \arg z < \frac{2(j+1)\pi }{2^n}\,\Big\}
\end{displaymath}
be the \emph{Luecking sets}.\par
\smallskip

The result is:

\begin{theoreme}[Luecking \cite{Lue}]\label{theo Luecking}
For every $p > 0$,  the composition operator $C_\phi$, assuming 
condition~\eqref{condition phi plus petit que un au bord}, is in the Schatten class $S_p$ if and only if:
\begin{equation}\label{critere Luecking}
\sum_{n\geq 0} 2^{np/2} \bigg( \sum_{j=0}^{2^n-1} \big[m_\phi (R_{n,j})\big]^{p/2} \bigg) <+\infty.
\end{equation}
\end{theoreme}

Majorizing $m_\phi (R_{n,j})$ by  
$m_\phi \big( W (\e^{2^{-n}(2j+1)i \pi}, 2^{-n})\big) \leq \rho_\phi (2^{-n})$, one gets:

\begin{corollaire}\label{coro Luecking}
Let $\phi \colon \D \to \D$ be an analytic self-map, with 
condition~\eqref{condition phi plus petit que un au bord}, and assume that $m_\phi$ is an $\alpha$-Carleson 
measure, with $\alpha > 1$. Then $C_\phi \in S_p$ for every $p > \frac{2}{\alpha - 1}\,\cdot$
\end{corollaire}

\noindent{\bf Proof.} Since $\rho_\phi (h) \lesssim h^\alpha$, one gets:
\begin{displaymath}
\sum_{n\geq 0} 2^{np/2} \bigg( \sum_{j=0}^{2^n-1} \big[m_\phi (R_{n,j})\big]^{p/2} \bigg) 
\lesssim \sum_{n \geq 0} 2^{np/2}.2^n.(2^{-n\alpha})^{p/2} 
= \sum_{n \geq 0} 2^{n [ 1 - (\alpha - 1)p/2]},
\end{displaymath}
which is $< +\infty$ since $1 - (\alpha - 1) \frac{p}{2} <0$. \qed
\bigskip

In order to get this corollary, we majorized crudely $m_\phi (R_{n,j})$ by $\rho_{\phi} (2^{-n})$. Actually, 
we shall see, through the results of this paper, that we lose too much with this majorization. Nevertheless, if we 
replace the Luecking sets by the \emph{dyadic Carleson windows}:
\begin{displaymath}
W_{n,j}=\Big\{ z\in\D \,;\  1 - 2^{-n} \leq |z| < 1\, ,\quad 
\frac{2j\pi}{2^n} \leq \arg(z) < \frac{2 (j + 1) \pi}{2^n} \Big\},
\end{displaymath}
($j=0, 1, \ldots, 2^n - 1$, $n=1,2,\ldots$), we have the same behaviour:

\begin{proposition}\label{Luecking vs Carleson}
Let $\mu$ be a finite positive measure on the open unit disk $\D$ and let $\alpha >0$. Then the following 
assertions are equivalent:\par
\smallskip
\qquad $(a)$ \quad $\displaystyle{
\sum_{n=1}^\infty \sum_{j=0}^{2^n - 1} 2^{n\alpha} \bigl(\mu(R_{n,j})\bigr)^\alpha < +\infty} $;\par
\smallskip
\qquad $(b)$ \quad $\displaystyle{
\sum_{n=1}^\infty \sum_{j=0}^{2^n - 1} 2^{n\alpha} \bigl(\mu(W_{n,j})\bigr)^\alpha < +\infty} $.
\end{proposition}

\noindent{\bf Proof.} It is clear that $(b)$ implies $(a)$ since $R_{n,j}\subset W_{n,j}$ for all $n$ and $j$ 
(we already used this in the proof of Corollary~\ref{coro Luecking}). \par
For the proof of the converse implication, we shall need the following sets, for positive integers $l, n$ 
with $l \geq n$, and $j\in \{0 ,1,\ldots, 2^n - 1\}$:
\begin{displaymath}
H_{l, n, j} = \Bigl\{ k\in \{0, 1, \ldots, 2^l -1\}\,;\ 
\frac{j}{2^n} \leq \frac{k}{2^l} < \frac{j + 1}{2^n} \Bigr\}\,.
\end{displaymath}

It is clear that we have, for every $n$ and $j$:
\begin{displaymath}
W_{n,j}  = \bigcup_{l \geq n}\; \bigcup_{k\in H_{l, n, j}} R_{l, k},
\end{displaymath}
and
\begin{displaymath}
\mu (W_{n, j}) = \sum_{l \geq n}\; \sum_{k \in H_{l, n, j}} \mu (R_{l, k})\,.
\end{displaymath}
\medskip

We shall first treat the case $\alpha\leq 1$, where we can use, for $x_1, x_2, \ldots, x_N \geq 0$:
\begin{displaymath}
(x_1 + x_2 + \cdots + x_N)^\alpha \leq x_1^\alpha + x_2^\alpha + \cdots + x_N^\alpha\,.
\end{displaymath}

We have:
\begin{align*}
\sum_{n=1}^\infty \sum_{j = 0}^{2^n - 1} 2^{n\alpha} \bigl(\mu(W_{n,j})\bigr)^\alpha
& = \sum_{n=1}^\infty \sum_{j=0}^{2^n - 1} 2^{n\alpha} \bigg( \sum_{l \ge n}\; 
\sum_{k \in H_{l,n,j}} \mu (R_{l,k}) \bigg)^\alpha \\
& \leq  \sum_{n=1}^\infty \sum_{j=0}^{2^n - 1} 2^{n\alpha} \sum_{l \geq n}\; 
\sum_{k\in H_{l,n,j}} \bigl(\mu(R_{l,k})\bigr)^\alpha \\
& =   \sum_{l =1}^\infty \sum_{k=0}^{2^l - 1} \;
\bigl(\mu(R_{l,k})\bigr)^\alpha  
\sum_{{(n,j) \,:\, n\le l \,;\, k\in H_{l,n,j}}} 2^{n\alpha}  
\end{align*}

Observe that, for every $n\le l$, there is only one $j$ such that $k\in  H_{l,n,j}$. Since we have
\begin{displaymath}
\sum_{n = 1}^{l} 2^{n\alpha} \le C_\alpha 2^{l\alpha},
\end{displaymath}
we get:
\begin{displaymath}
\sum_{n=1}^\infty \sum_{j=0}^{2^n - 1} 2^{n\alpha} \bigl(\mu(W_{n,j})\bigr)^\alpha
\le C_\alpha
\sum_{l=1}^\infty \sum_{k=0}^{2^l - 1} 2^{l\alpha} \bigl(\mu(R_{l,k})\bigr)^\alpha\,,
\end{displaymath}
and $(a)$ implies $(b)$ in the case $\alpha\le 1$.
\medskip

If $\alpha>1$, we can use H\"older's inequality. Let $\beta$ be the conjugate exponent of $\alpha$. Choose 
$a$ such that $1 < a < 2 < a^\beta$. Then:
\begin{align*}
\mu (W_{n,j}) 
& = \sum_{l\ge n \,:\, k\in H_{l,n,j}}\; \mu (R_{l,k}) \\
& \le \bigg( \sum_{l \ge n \,:\, k\in H_{l,n,j}}\; a^{-l\beta} \bigg)^{1/\beta}
\bigg( \sum_{l \ge n \,:\, k\in H_{l,n,j}}\; a^{l\alpha} 
\bigl(\mu (R_{l,k})\bigr)^\alpha \bigg)^{1/\alpha}
\end{align*}
Using that $|H_{l,n,j}|=2^{l-n}$, we get:
\begin{displaymath}
\bigg[\sum_{l\ge n \,:\, k\in H_{l,n,j}}\; a^{-l\beta}\bigg]^{1/\beta} = 
\bigg[\sum_{l \ge n} 2^{l - n} a^{- l\beta}\bigg]^{1/\beta} = 
\biggl( 2^{-n} \frac{(2 a^{-\beta})^n}{1 - 2 a^{-\beta}}\biggr)^{1/\beta}= C_\beta a^{-n} .
\end{displaymath}
Therefore we have:
\begin{align*}
\sum_{n=1}^\infty \sum_{j=0}^{2^n - 1} 2^{n\alpha} \bigl(\mu(W_{n,j})\bigr)^\alpha
&\le \sum_{n=1}^\infty \sum_{j=0}^{2^n - 1} 2^{n\alpha} C_\beta^\alpha a^{-n\alpha}        
\sum_{l \ge n \,:\, k\in H_{l,n,j}}\; a^{l\alpha} \bigl(\mu (R_{l,k})\bigr)^\alpha \\
& = C_\beta^\alpha \sum_{l =1}^\infty \sum_{k=0}^{2^l - 1} \;
\bigl(\mu(R_{l,k})\bigr)^\alpha a^{l \alpha} \sum_{{(n,j) \,:\, n\le l\,, k\in H_{l,n,j}}} (2/ a) ^{n\alpha}  \cr
& \lesssim \sum_{l =1}^\infty \sum_{k=0}^{2^l - 1} \;
\bigl(\mu(R_{l,k})\bigr)^\alpha a^{l \alpha} (2/a)^{l \alpha} \\
& = \sum_{l=1}^\infty \sum_{k=0}^{2^l - 1} 2^{l\alpha} \bigl(\mu (R_{l,k})\bigr)^\alpha\,.
\end{align*}
We have hence proved that $(a)$ implies $(b)$ for $\alpha > 1$ and therefore 
Proposition~\ref{Luecking vs Carleson} follows. \qed
\medskip

As a corollary we prove a necessary condition that $\rho_\phi$ must satisfy when $C_\phi$ is in the Schatten 
class $S_p$.

\begin{proposition}\label{condition necessaire Schatten}
Let $\phi \colon\D\to \D$ be an analytic self-map. If the composition operator
$C_\phi \colon H^2 \to H^2$ is in the Schatten class $S_p$ for some $p> 0$, then, as $h$ goes to $0$:
\begin{equation}\label{limite condition necessaire Schatten}
\rho_\phi (h) = o\, \bigg( h \Big( \log \frac{1}{h} \Big)^{-2/p} \bigg).
%\lim_{h \to 0^+} \Big( \frac{\rho_\phi (h)}{h}\Big)^{p/2} \log (1/h) = 0.
\end{equation}
\end{proposition}

\noindent{\bf Proof.} Thanks to Luecking's characterization and the equivalence in 
Proposition~\ref{Luecking vs Carleson}, we have, for the pullback mesure $m_\phi$:
\begin{equation}\label{eq 1 condition necessaire}
\sum_{n=1}^\infty \sum_{j=0}^{2^n - 1} 2^{n p/2} \big(m_\phi (W_{n,j}) \big)^{p/2} < +\infty.
\end{equation}
Observe that, for $h=2^{-n}$ each window $W (\xi, h)$ is contained in the union of at most three of the 
$W_{n, j}$'s; hence:
\begin{displaymath}
\big(\rho_\phi (2^{-n})\big)^{p/2} \leq 
\Big(3 \max_{0\le j\le 2^n - 1} m_\phi (W_{n,j}) \Big)^{p/2} 
\le 3^{p/2} \sum_{j=0}^{2^n - 1} \big(m_\phi (W_{n,j}) \big)^{p/2}  \,,
\end{displaymath}
and \eqref{eq 1 condition necessaire} yields:
\begin{displaymath}
\sum_{n=1}^\infty \bigl(\rho_\phi (2^{-n})\bigr)^{p/2} 2^{n p/2} < +\infty\,.
\end{displaymath}
Hence, setting:
\begin{equation}
\gamma_n = \sum_{n/2\le k \le n} \bigl(\rho_\phi (2^{-k})\bigr)^{p/2} 2^{kp/2}\,,
\end{equation}
we have:
\begin{equation}\label{eq 2 condition necessaire}
\lim_{n\to\infty} \gamma_n = 0\,.
\end{equation}
\medskip

Now, by using \cite{CompOrli}, Theorem~4.19, we get a constant $C > 0$ such that, for $k\le n$:
\begin{displaymath}
C \rho_\phi (2^{-k}) \ge 2^{n - k} \rho_\phi (2^{-n})\,,
\end{displaymath}
and so:
\begin{equation}\label{eq 3 condition necessaire}
C^{p/2} \gamma_n \ge (n/2) \Big(\frac{\rho_\phi(2^{-n})}{2^{-n}}\Big)^{p/2} \,.
\end{equation}
\par

To finish the proof, it remains to consider, for every $h\in (0,1/2)$, the integer $n$ such that 
$2^{-n-1} < h \le 2^{- n}$; then \eqref{eq 2 condition necessaire} and \eqref{eq 3 condition necessaire} give:
\begin{displaymath}
\lim_{h\to 0^+} \Big( \frac{\rho_\phi (h)}{h}\Big)^{p/2} \log(1/h) = 0,
\end{displaymath}
as announced. \qed
\medskip

\noindent{\bf Remark.} We can also deduce Corollary~\ref{coro Luecking} from the following result.

\begin{proposition}
If $\mu$ is a $\beta$-Carleson probability measure on $\D$, with $\beta > 2$, then the Poisson integral 
${\cal P} \colon L^2 (\T) \to L^2 (\mu)$ is in the Schatten class $S_p$ for any $p > 2/ (\beta - 1)$.
\end{proposition}

\noindent{\bf Proof.} We may assume that $p \leq 2$ since $S_{p_1} \subseteq S_{p_2}$ when 
$p_1 \leq p_2$. To have ${\cal P} \in S_p$, it suffices then to have 
$\sum_{n\in \Z} \| {\cal P}(e_n) \|_{L^2 (\mu)}^p < +\infty$, where $e_n (\e^{it}) = \e^{int}$ 
(see \cite{Koenig}, Proposition~1.b.16, page 40, for example).\par
But $({\cal P} e_n ) (z) = z^n$ for $n \geq 0$ and $({\cal P} e_n ) (z) = \bar{z}^{|n|}$ for $n \leq - 1$. Hence:
\begin{displaymath}
\sum_{n\in \Z} \| {\cal P}(e_n) \|_{L^2 (\mu)}^p
\leq 2\sum_{n = 0}^\infty \Big( \int_{\D} |z|^{2n}\,d\mu \Big)^{p/2} .
\end{displaymath}
But
\begin{align*}
\int_{\D} |z|^{2n}\,d\mu 
& = \int_0^1 2n t^{2n -1} \mu (|z | \geq t) \,dt \\
& = 2 n \int_0^1 (1 - x)^{2n - 1} \mu (|z | \geq 1 - x) \,dx \\
& \lesssim n \int_0^1 (1 - x)^{2n - 1} \frac{1}{x}\, x^\beta \,dx \\
& \quad \text{since $\{ | z | \geq 1 - x \}$ can be split in $O (1/x)$ windows $W (a, x)$} \\
& = n \int_0^1 (1 - x)^{2n - 1} x^{\beta - 1}\,dx  \lesssim n.n^{-\beta}.
\end{align*}
Hence
\begin{displaymath}
\sum_{n\in \Z} \| {\cal P}(e_n) \|_{L^2 (\mu)}^p \lesssim \sum_{n \geq 1} \frac{1}{n^{(\beta - 1)p/2}}
\raise 1,5pt \hbox{,}
\end{displaymath}
which is finite since $(\beta - 1)p/2 > 1$. \qed

\bigskip

%%%%%%%%%%%%%%%%%%%%%%%%%

\subsection{A general contruction}\label{general construction}

In this subsection, we are going to describe a general way to construct symbols with some prescribed conditions. 
A particular case of this construction has been used in \cite{CompOrli}, Theorem~4.1. We also shall use it in 
\cite{Orlicz-Hardy}.
\medskip

Let
\begin{displaymath}
f (t) = \sum_{k=0}^\infty a_k \cos (kt)
\end{displaymath}
be an even, non-negative, $2\pi$-periodic continuous function, vanishing at the origin: $f (0) = 0$, and such 
that:
\begin{equation}\label{deuxieme condition}
\text{$f$ is \emph{strictly} increasing on $[0,\pi]$.}
\end{equation}

The Hilbert transform (or conjugate function) ${\cal H}f$ of $f$ is:
\begin{displaymath}
{\cal H}f (t) =\sum_{k=1}^\infty a_k \sin (kt).
\end{displaymath}
We shall assume moreover that, as $t$ tends to zero:
\begin{equation}\label{derivee de Hf}
\big({\cal H}f\big)' (t) = o\, ( 1/t^2).
\end{equation}
\par

Let now $F \colon \D \to \Pi^+ = \{ \Re z >0\}$ be the analytic function whose boundary values are
\begin{equation}\label{def grand F}
F^\ast (\e^{it}) = f (t) + i {\cal H}f (t).
\end{equation}
One has:
\begin{displaymath}
\qquad F (z) = \sum_{k=0}^\infty a_k z^k, \qquad |z| < 1.
\end{displaymath}
\par
We define:
\begin{equation}\label{def grand Phi}
\qquad \Phi (z) = \exp \big( - F(z)\big), \qquad z\in \D.
\end{equation}
Since $f$ is non-negative, one has
\begin{displaymath}
\Re F (z) = \frac{1}{2\pi} \int_{-\pi}^{\pi} f (t) P_z (t)\,dt > 0,
\end{displaymath}
so that $| \Phi (z) | < 1$ for every $z \in \D$: $\Phi$ is an analytic self-map of $\D$, and 
$|\Phi^\ast | = \exp (- f ) < 1$ \emph{a.e.} . Note that the assumption $f (0) = 0$ means that 
$\Phi^\ast (1) = 1$; we then have $\| \Phi \|_\infty = 1$, which is necessary for the non-compactness of 
$C_\Phi$.\par
\smallskip

An example is $\Phi (z ) = \exp \big( - (1 - z)/2\big)$, for which $f (t) = \sin^2 (t/2)$.\par
\smallskip

\begin{lemme}\label{premier lemme}
Assume that $f$ and ${\cal H}f$ are ${\cal C}^1$ functions.
Then, the Carleson function $\rho_\Phi$ of $\Phi$ is not $o\,(h)$ when $h$ goes to $0$, and so the 
composition operator $C_\Phi \colon H^2 \to H^2$ is not compact.
\end{lemme}

Note that the hypothesis of the lemma holds, for example, when:
\begin{equation}\label{premiere condition}
\sum_{k=0}^\infty k\,|a_k| < +\infty.
\end{equation}

\noindent{\bf Proof.} The non-compactness of $C_\Phi$ follows immediately from the ``angular 
derivative'' condition (\cite{Shap-livre}, Theorem~3.5: see \eqref{angular derivative} in the remark at the 
end of the next section). But we shall give a proof using the Carleson function, in order to illustrate the 
methods to be used later on.\par
It will be enough to minorize $\mu_\Phi [W (1,h)]$. Since $f$ and ${\cal H}f$ are 
${\cal C}^1$, we have $| f (t) | \leq C |t|$ and $| {\cal H}f (t) | \leq C |t|$ for some positive constant $C$. 
Now, if $| t | \leq h/ C$, we see that 
\begin{displaymath}
| \Phi^\ast ( \e^{it}) | = \e^{ - f (t) } \geq \e^{- C |t|} \geq \e^{-h} \geq 1 - h,
\end{displaymath}
and $| \arg \Phi^\ast (\e^{it}) | = | {\cal H} f (t) | \leq C |t| \leq h$; hence $\Phi^\ast (\e^{it}) \in W (1,h)$. 
This means that
\begin{displaymath}
m_\Phi [W (1, h)] \geq m (\{ \e^{it}\,;\ | t | \leq h / C \}) \geq \frac{1}{\pi C}\, h,
\end{displaymath}
and the lemma follows. \qed
\bigskip

We shall now perturb $\Phi$ by considering
\begin{equation}\label{M}
\qquad M (z) = \exp\Big( - \frac{1+ z }{ 1 -z }\Big)\,, \qquad |z| <1,
\end{equation}
and 
\begin{equation}\label{def general phi}
\qquad \phi ( z) = M (z) \Phi (z), \qquad | z | <1.
\end{equation}
One has:
\begin{displaymath}
\phi^\ast (\e^{it})  = \e^{ - f(t)}\,\e^{- i ({\cal H} f (t) + \cot \frac{t}{2} )}.
\end{displaymath}

We will now, according to the various choices of $f$, study the behaviour of $\phi$ with respect to the 
Carleson windows.\par
\smallskip

We set:
\begin{equation}
\gamma (t ) =  {\cal H} f (t) + \cot \frac{t}{2} = \frac{2}{t} + r (t),
\end{equation}
where the derivative of the odd function $r$ satisfies $r ' (t) = o\, (1/t^2)$ .\par
Now, we have:

\begin{lemme}\label{deuxieme lemme}
When $h > 0$ goes to $0$, one has:
\begin{equation}
h f^{-1} (h) \lesssim \rho_\phi (h) \lesssim h f^{-1} (2h).
\end{equation}
\end{lemme}

\noindent{\bf Proof.} Let $a = \e^{i\theta} \in \T$, $|\theta | \leq \pi$. We may assume that 
$0 < h \leq h_0 \leq 1/2$, for some $h_0$ small enough. \par
We have to analyze the set of $t$'s such that 
$\phi^\ast (\e^{-it}) = \e^{- f (t)} \e^{i \gamma (t)} \in W (a, h)$, which 
imposes two constraints. Without loss of generality, we may analyze only the set of positive $t$'s, \emph{i.e.} 
$0 < t \leq \pi$.\par
\smallskip

\emph{Modulus constraint.} We must have $| \phi^\ast (\e^{-it}) | \geq 1 - h$, \emph{i.e.} 
$\e^{- f (t)} \geq 1 - h$, or $f (t) \leq \log \frac{1}{1 - h}$, which is $\leq 2h$ since $h \leq 1/2$. Hence we 
must have:
\begin{equation}\label{contrainte module}
0 < t \leq f^{-1} (2h).
\end{equation}

\emph{Argument constraint.} We must have $| \gamma (t) - \theta | \leq h \mod 2\pi$, \emph{i.e.}, since 
we have assumed that $t >0$:
\begin{equation}\label{intervalle}
\gamma (t) \in \bigcup_{n \geq 0} [ \theta - h + 2n \pi, \theta + h + 2n \pi] = \bigcup_{n \geq 0} J_n (h).
\end{equation}
Since $\gamma (t) \to +\infty$ as $t \mathop{\to}\limits^{\scriptscriptstyle >} 0$, and since we 
already have $t \leq f^{-1}(2h)$ by \eqref{contrainte module}, we know that $\gamma (t) > 2 \pi$ for $h$ 
small enough; hence we have $\gamma (t) \in J_n (h)$ only for $n \geq N_h$, where the integer $N_h$ 
goes to infinity when $h$ goes to $0$; in particular, we may assume that, for $h$ small enough, we have 
$\gamma (t) \in J_n (h)$ for $n \geq 1$ only.\par
Let $I_n (h) = \gamma^{-1} \big( J_n (h) \big)$. Since $\gamma (t) = \frac{2}{t} + r (t)$, and 
$r' (t) = o\, (1/t^2)$, $\gamma (t)$ is decreasing for $0 < t \leq h_0$, for $h_0$ small enough. Hence, for 
$h$ small enough:
\begin{displaymath}
I_n (h)  = \big[ \gamma^{-1} (\theta +h + 2n \pi), \gamma^{-1} (\theta - h + 2n \pi) ].
\end{displaymath}
Since $\gamma (t) = \frac{2}{t} + o\, (1/t)$ and $|\gamma ' (t) | = \frac{2}{t^2} + o\, (1/t^2)$, one has:
\begin{equation}\label{encadrement intervalle}
\frac{c_1}{n} \leq \min I_n (h) \leq \max I_n (h) \leq \frac{c_2}{n}\,\raise 1,5pt \hbox{,}
\end{equation}
where $c_1$ and $c_2$ are two universal positive constants. By the mean-value theorem, we get that:
\begin{align*}
2 \pi m \big( I_n (h)\big) 
& = \gamma^{-1} (\theta - h + 2n \pi) - \gamma^{-1} (\theta +h + 2n \pi) \\
& = 2h | (\gamma^{-1})' (\xi_n)| = \frac{2h }{|\gamma ' (t_n)|}\,\raise 1,5pt \hbox{,}
\end{align*}
for some $\xi_n \in J_n (h)$ and $t_n \in I_n (h)$. But, \eqref{encadrement intervalle} ensures that 
$\frac{c_1}{n} \leq t_n \leq \frac{c_2}{n}$ and, since $|\gamma ' (t) | = \frac{2}{t^2} + o (1/t^2)$, we get 
that:
\begin{equation}\label{mesure intervalle}
m\big( I_n (h)\big) \approx \frac{h}{n^2}\cdot
\end{equation}
\smallskip

Now:\par

1) Assume that $\phi^\ast (\e^{-it}) \in W (a, h)$. By \eqref{intervalle}, \eqref{encadrement intervalle} and 
\eqref{contrainte module}, we must have $t \in I_n (h)$ with $\frac{c_1}{n} \leq t \leq f^{-1}(2h)$. Hence, 
if $n_0$ is the integer part of $\frac{c_1}{f^{-1}(2h)}$, we must have $n \geq n_0$. Now, 
\eqref{mesure intervalle} shows that:
\begin{align*}
m\big( \{ t\in ]0,\pi]\,;\ \phi^\ast (\e^{- it}) \in W (a,h)\}\big) 
& \leq \sum_{n \geq n_0} m\big( I_n (h)\big) \\
& \lesssim \sum_{n \geq n_0} \frac{h}{n^2} \lesssim \frac{h}{n_0} \lesssim h f^{-1}(2h).
\end{align*}
\par

2) We want to minorize, as it suffices, $m\big( \{ t\in ]0,\pi]\,;\ \phi^\ast (\e^{it}) \in W (1,h)\}\big)$. Let 
$n_1$ be the integer part of $\frac{c_2}{f^{-1}(h)} +1$; we have:
\begin{equation}\label{implication fenetre}
t \in \bigcup_{n \geq n_1} I_n (h) \qquad \Longrightarrow \qquad  \phi^\ast (\e^{it}) \in W (1,h)
\end{equation}
because $t \in I_n (h)$ for $n \geq n_1$ implies $t \leq \frac{c_2}{n_1} \leq f^{-1} (h)$, so that 
$| \phi^\ast (\e^{it}) | = \e^{ - f(t)} \geq \e^{-h} \geq 1 - h$, and the modulus constraint for 
$\phi^\ast (\e^{it}) $ is automatically satisfied. Since the argument constraint is satisfied by construction, as 
$t$ belongs to some $I_n(h)$, this proves \eqref{implication fenetre}.\par
As a consequence, we have, using \eqref{mesure intervalle}:
\begin{displaymath}
\rho_\phi (h) \geq m_\phi \big( W (1, h) \big) \geq \sum_{n \geq n_1} m\big( I_n (h)\big) \gtrsim 
\sum_{n \geq n_1} \frac{h}{n^2} \gtrsim \frac{h}{n_1} \gtrsim h f^{-1} (h),
\end{displaymath}
and this ends the proof of Lemma~\ref{deuxieme lemme}. \qed
\bigskip

%%%%%%%%%%%%%%%%%%%%%%%%%%%%%

\section{Composition operators with symbol of same modulus}\label{section same modulus}

J. Shapiro and P. Taylor (\cite{Shap-Tay}, Theorem~3.1) characterized Hilbert-Schmidt composition 
operators $C_\phi$ on $H^2$ (\emph{i.e.} $C_\phi\in S_2$), and this characterization depends only on 
the modulus of $\phi^\ast$ on $\partial \D$. It follows that if $\phi_1$ and $\phi_2$ are two symbols such 
that $|\phi_1^\ast | = |\phi_2^\ast|$, and $C_{\phi_2}$ is Hilbert-Schmidt, then $C_{\phi_1}$ is also 
Hilbert-Schmidt, and in particular, compact.\par
It appears that this is a limiting case, as we shall see in Theorem~\ref{meme module Schatten}. 
\medskip

Actually, in view of Corollary~\ref{coro Luecking}, we may formulate this problem in the following way. 
Assume that the composition operator $C_{\phi_1}$ is not compact on $H^2$, and that 
$|\phi_1^\ast | = |\phi_2^\ast|$; for which values of $\alpha$, can $m_{\phi_2}$ be an 
$\alpha$-Carleson measure?\par
Note that if $m_{\phi_2}$ is an $\alpha$-Carleson measure, we necessarily must have 
$\alpha \leq 2$. Indeed, if $\alpha >  2$, Corollary~\ref{coro Luecking} implies that $C_{\phi_2}\in S_2$, 
and hence $C_{\phi_1} \in S_2$ as well.\par
\medskip

When $\alpha < 2$, the situation is very different, and one has:

\begin{theoreme}\label{meme module puissance}
For every $\alpha$ with $1 < \alpha < 2$, there exist two symbols $\phi_1$ and $\phi_2$ having the 
same modulus on $\partial \D$ and such that $\rho_{\phi_1}(h) \approx h$, but 
$\rho_{\phi_2}(h) \approx h^\alpha$.
\end{theoreme}

It follows from Corollary~\ref{coro Luecking} that:

\begin{theoreme}\label{meme module Schatten}
For every $p > 2$, there exist two symbols $\phi_1$ and $\phi_2$ having the same modulus on $\partial \D$ 
and such that $C_{\phi_1}$ is not compact on $H^2$, but $C_{\phi_2}$ is in the Schatten class $S_p$. 
\end{theoreme}

\noindent{\bf Proof of Theorem~\ref{meme module puissance}.}  In \cite{CompOrli}, Theorem~4.1, we proved 
a particular case of this result, corresponding to $\alpha =3/2$. We took there 
$\phi_1 (z) =\frac{1+z}{2}\,$. This function behaves as 
$\exp\big( - [\sin^2 (t/2)] + i {\cal H} [\sin^2 (t/2)]\big)$, 
and to prove Theorem~\ref{meme module puissance}, we shall just change the power $2$ of $\sin (t/2)$.
\par\smallskip

We shall use the following lemma, whose proof will be postponed.

\begin{lemme}\label{lemme coeff}
For $0 < \beta < 2$, let $f (t) = |\sin \frac{t}{2}|^\beta$. Then:
\begin{displaymath}
f (t) = c_0 + \sum_{k=1}^\infty c_k \cos kt\,,
\end{displaymath}
with:
\begin{equation}\label{signe coeff}
c_k < 0 \quad \text{for all } k\geq 1;
\end{equation}
and
\begin{equation}\label{estim coeff}
c_k = O\,\Big(\frac{1}{k^{\beta +1}}\Big)\,\cdot
\end{equation}
In particular, for $\beta >1$, the series $\sum_{k\geq 1} k c_k$ is convergent, and 
$\sum_{k\geq 1} k c_k  < 0$.
\end{lemme}
\par\medskip

Taking $\beta = \frac{1}{\alpha - 1}$\,, which is $>1$, and $f (t) = | \sin (t/2)|^\beta$ as in 
Lemma~\ref{lemme coeff}, observe that $f$ satisfies the assumptions \eqref{premiere condition} and 
\eqref{deuxieme condition} of the general construction of the Subsection~\ref{general construction} (note that 
for $\beta \geq 2$, $f$ is a ${\cal C}^2$ function, and hence, we have \eqref{premiere condition} directly, 
without using Lemma~\ref{lemme coeff}). With the notation of that subsection (see \eqref{def grand F}, 
\eqref{def grand Phi}, and \eqref{M}), set:

\begin{displaymath}
\phi_1 =  \Phi \qquad \text{and} \qquad \phi_2 = M \Phi.
\end{displaymath}
One has:

\begin{displaymath}
|\phi_1^\ast| = | \phi_2^\ast | \quad a.e.
\end{displaymath}
\par

Lemma~\ref{premier lemme} and MacCluer's theorem show that $C_{\phi_1}$ is not compact on $H^2$. 
On the other hand, since $f^{-1}(h) \approx h^{1/\beta}$, Lemma~\ref{deuxieme lemme} shows that 
$\rho_{\phi_2} (h) \approx h^\alpha$, with $\alpha = 1 + \frac{1}{\beta} \in ]1,2[$.\par
\smallskip

That ends the proof of Theorem~\ref{meme module puissance}. \qed
\bigskip

\noindent{\bf Proof of Lemma~\ref{lemme coeff}.} Before beginning the proof, it should be remarked that for 
$\beta =2$ (\emph{i.e.} $\alpha =3/2$, which is the case processed in \cite{CompOrli}, Theorem~4.1), we 
have a trivial situation: $\sin^2 \frac{t}{2} = \frac{1}{2} - \frac{1}{2} \cos t$. \par
\emph{Proof of \eqref{signe coeff}}. For $0 < p < 1$, we have the well-known binomial expansion, for 
$-1 \leq x \leq 1$:
\begin{displaymath}
(1 -x)^p = 1 - \sum_{k=1}^{\infty} \alpha_k x^k,
\end{displaymath}
with:
\begin{displaymath}
\alpha_k = \frac{p (1 - p) \cdots (k - 1 -p)}{k!} >0.
\end{displaymath}
Taking $x = \cos t$ and $p= \beta /2$, we get:
\begin{equation}\label{etoile}
2^{\beta/2} \Big| \sin {\frac{t}{2}} \Big|^\beta = 1 - \sum_{k=1}^\infty \alpha_k (\cos t)^k.
\end{equation}
Now, we know that:
\begin{displaymath}
(\cos t)^k = \sum_{j=0}^k b_{j,k} \cos (k - 2j) t \,\raise 1,5 pt \hbox{,}
\end{displaymath}
with $b_{j,k} > 0$. Substituting this in \eqref{etoile}, grouping terms, and dividing by $2^{\beta/2}$, we 
get \eqref{signe coeff}.\par

\emph{Proof of \eqref{estim coeff}}. We shall separate the cases $\beta = 1$, $0 < \beta < 1$, and 
$1 < \beta < 2$.\par
For $\beta =1$, one has, in an explicit way:

\begin{displaymath}
\Big|\sin \frac{t}{2} \Big| =\frac{2}{\pi} - \frac{4}{\pi} \sum_{k=1}^\infty \frac{\cos kt}{4 k^2 - 1} 
=c_0 +\sum_{k=1}^\infty c_k \cos kt \,.
\end{displaymath}
\smallskip

Assume $0 < \beta < 1$. Since $\frac{\sin t/2}{t/2} >0$ on $[0,\pi]$, we can write:

\begin{displaymath}
\Big( \frac{\sin t/2}{t/2}\Big)^\beta = 1  + t^2 u(t)\,,
\end{displaymath}
where $u$ is a ${\cal C}^\infty$ function on  $[0,\pi]$. Then:

\begin{displaymath}
c_k = \frac{2}{2^\beta \pi} \int_0^\pi t^\beta \cos kt \,dt + 
\frac{2}{2^\beta \pi} \int_0^\pi t^{\beta +2} u(t) \cos kt \,dt.
\end{displaymath}
The second integral is $O (k^{-2})$, as easily seen by making two integrations by parts. The first one 
writes:

\begin{align*}
\int_0^\pi t^\beta \cos kt \,dt 
& = -\frac{\beta}{k} \int_0^\pi t^{\beta - 1} \sin kt \,dt \\
& =  -\frac{\beta}{k} \int_0^{k\pi} \Big(\frac{x}{k}\Big)^{\beta - 1} \sin x \,\frac{dx}{k} \\
& \sim -\frac{\beta}{k^{\beta+1}} \int_0^{+\infty} x^{\beta - 1} \sin x\,dx.
\end{align*}
This last integral is convergent and positive. Hence, since $\beta + 1 < 2$:
\begin{displaymath}
c_k \sim - \delta k^{-(\beta + 1)},
\end{displaymath}
where $\delta$ is a positive constant.
\par

Before continuing, let us observe that we have similarly, due to the vanishing of the integrated terms:

\begin{equation}\label{grand O integrale}
\int_0^\pi \Big(\sin \frac{t}{2}\Big)^\sigma \sin (2k+1) \frac{t}{2}\,dt 
= O (k^{-(\sigma + 1)}) + O (k^{-2}) = O (k^{-(\sigma + 1)}), 
\end{equation}
for $0 < \sigma < 1$.
\smallskip

Assume now  $1 < \beta < 2$. We have:

\begin{align*}
\frac{\pi}{2} c_k 
& = \int_0^\pi f (t) \cos kt\,dt = - \frac{1}{k} \int_0^\pi f '(t) \sin kt\,dt \\
& = - \frac{\beta}{2k} \int_0^\pi \big(\sin (t/2)\big)^{\beta - 1}  \cos (t/2) \sin (kt)\,dt \\
& = - \frac{\beta}{4k} \bigg[\int_0^\pi \Big(\sin \frac{t}{2} \Big)^{\beta - 1}\sin (2k+1)\frac{t}{2} \,dt \\
& \qquad \qquad \qquad \qquad \qquad 
+ \int_0^\pi \Big(\sin \frac{t}{2} \Big)^{\beta - 1} \sin (2k - 1) \frac{t}{2} \,dt \bigg] \\
& = \frac{1}{k}\, O (k^{-\beta}) = O (k^{-(\beta +1)}), 
\end{align*}
in view of \ref{grand O integrale}, applied with $\sigma = \beta - 1 \in ]0,1[$.\par
This ends the proof of Lemma~\ref{lemme coeff}. \qed
\bigskip

\noindent{\bf Remark.} The following question arises naturally: does $C_\phi \in S_p$ for some $p < \infty$ 
imply that $\mu_\phi$ is $\alpha$-Carleson for some $\alpha >1$? We shall see in the next section 
(Remark~2 after the proof of Proposition~\ref{prop Schatten}) that the answer is negative.\par
\smallskip

Another question, related to our work, has been raised by K. Kellay: given a compact composition operator 
$C_\phi \colon H^2 \to H^2$, and another symbol $\psi \colon \D \to \D$, is the composition 
operator $C_{\phi \psi} \colon H^2 \to H^2$ still compact? This is the case if $\phi \psi$ is univalent. 
In fact, the compactness of $C_\phi$ implies (\cite{Shap-livre}, Theorem~3.5) that:
\begin{equation}\label{angular derivative}
\lim_{|z| \to 1} \frac{1 - |\phi (z)|}{1 -|z|} = +\infty; 
\end{equation}
hence, since $|\phi (z) \psi (z) | \leq |\phi (z)| $ for every $z \in \D$, we have:
\begin{displaymath}
\lim_{|z| \to 1} \frac{1 - |\phi (z) \psi (z)|}{1 -|z|} = +\infty,
\end{displaymath}
which implies the compactness of $C_{\phi\psi}$, thanks to the univalence of $\phi \psi$ 
(\cite{Shap-livre}, Theorem~3.2).\par

In \cite{CompOrli}, Proposition 4.2, we proved a related result: let $\phi_1, \phi_2 \colon \D \to \D$ be 
univalent analytic self-maps such that $|\phi_1^\ast| \leq |\phi_2^\ast|$ on $\partial\D$; if 
$C_{\phi_2} \colon H^2 \to H^2$ is compact, and $\phi_2$ vanishes at some point $a\in \D$, then 
$C_{\phi_1} \colon H^2 \to H^2$ is also compact. Note that the vanishing condition for $\phi_2$ 
in that result is automatic (since $\phi_2$ must have a fixed point $a\in \D$ because the 
composition operator $C_{\phi_2}$ is compact, see \cite{Shap-livre}, page~84, \S~5.5, Corollary), but the 
univalence condition for $\phi_1$ and $\phi_2$ cannot be dropped. In fact, take  
$\phi_1 (z) = \exp \big(- (1 - z)/2\big)$ and $\phi_2 (z) = z\, M (z)\,\exp \big(- (1 - z)/2\big)$ (where $M (z)$ 
is defined by \eqref{M}). One has $|\phi_1^\ast| = |\phi_2^\ast|$ on $\partial \D$, and $\phi_1$ is univalent 
(if $\phi_1 (z) = \phi_1 (w)$, $k = 0$ is the only integer such that $\frac{z}{2} - \frac{w}{2} = 2k\pi i$, since  
$|z| + |w| \leq 2$). But $\phi_1$ is the function $\Phi$ defined by \eqref{def grand Phi}, with 
$f (t) = \sin^2 (t/2)$. Since $f$ and ${\cal H}f$ are ${\cal C}^1$, Lemma~\ref{premier lemme} says that 
$C_{\phi_1}$ is not compact. On the other hand, one has $f^{-1} (h) \approx h^{1/2}$, so 
Lemma~\ref{deuxieme lemme} gives the compactness of $C_{\phi_2}$ (and even, $C_{\phi_2} \in S_p$ 
for every $p >4$, by Corollary~\ref{coro Luecking}).\par
It should be pointed out that, however, $\phi_1$ cannot be written $\phi_1 = \phi_2 \psi$ for some analytic 
self-map $\psi \colon \D \to \D$.

%%%%%%%%%%%%%%%%

\section{Composition operators in Schatten classes}\label{log puissance theta}

In \cite{Shap-Tay}, Theorem~4.2, J. Shapiro and P. Taylor constructed a family of composition operators 
$C_{\phi_\theta} \colon H^2 \to H^2$, indexed by a parameter $\theta >0$ such that  
$C_{\phi_\theta}$ is always compact, but $C_{\phi_\theta}$ is Hilbert-Schmidt if and only if 
$\theta > 2$. In this section, we shall slightly modify the symbol $\phi_\theta$, 
and shall study the membership of $C_{\phi_\theta}$ in the Schatten classes $S_p$. In \cite{Orlicz-Hardy}, 
we study on which Hardy-Orlicz spaces $H^\Psi$ these composition operators $C_{\phi_\theta}$ 
are compact.\par

\begin{theoreme}\label{theo Schatten} 
For every $p_0 > 0$, there exists an analytic self-map $\phi \colon \D \to\D$ such that the composition 
operator $C_\phi \colon H^2 \to H^2$ is in the Schatten class $S_p$ for every $p > p_0$, but not in 
$S_{p_0}$.
\end{theoreme}

\noindent{\bf Proof.} We shall use the same function as J. Shapiro and P. Taylor in 
\cite{Shap-Tay}, \S~4, with slight modifications. This modified function will be easier to analyze.\par
Let $\theta > 0$.\par
For $\Re z > 0$, $\log z$ will be the principal determination of the logarithm. Let, for $\eps >0$:
\begin{equation}\label{V indice epsilon}
V_\eps = \{ z\in \C\,;\ \Re z >0 \ \text{and}\ |z | < \eps\}.
\end{equation}
and consider, for $\eps >0$ small enough:
\begin{equation}\label{f indice theta}
f_\theta (z) = z (- \log z )^\theta, \qquad z \in V_\eps.
\end{equation}

\begin{lemme}\label{definition f indice theta}
For $\eps >0$ small enough, one has $\Re f_\theta (r \e^{i \alpha}) > 0$ for $0 < r < \eps$ and
$| \alpha | < \pi/2$. Moreover, one has $\Re f_\theta^\ast ( z) > 0$ for all 
$z \in \partial V_\eps \setminus\{0\}$.
\end{lemme}

\noindent{\bf Proof.} Actually, $f_\theta$ can be defined on $\overline{V_\eps} \setminus \{0\}$, and we 
shall do that.\par
Let, for $| \alpha | \leq \pi/2$ and $0 < r \leq \eps$:
\begin{displaymath}
Z_\alpha = \big(- \log (r \e^{i\alpha}) \big)^\theta = \Big( \log \frac{1}{r} - i \alpha \Big)^\theta.
\end{displaymath}
One has $f_\theta (r \e^{i\alpha}) = r \e^{i\alpha} |Z_\alpha| \e^{i \arg Z_\alpha}$, so that:
\begin{displaymath}
\Re f_\theta (r \e^{i\alpha}) = r\, |Z_\alpha| \cos (\alpha + \arg Z_\alpha).
\end{displaymath}
On the other hand,
\begin{displaymath}
\arg Z_\alpha = - \theta  \arctan \frac{\alpha}{\log 1/r}\, \cdot
\end{displaymath}
Since $\arctan x \geq x/2$ for $ 0 \leq x \leq 1$, we get, for $0 < r \leq \eps \leq \e^{-\pi/2}$:
\begin{align*}
| \alpha + \arg Z_\alpha | 
& \leq |\alpha| \Big( 1 -  \frac{\theta}{\alpha}  \arctan \frac{\alpha}{\log 1/r} \Big)  \\
& \leq \frac{\pi}{2} \Big( 1 - \frac{\theta}{2\log 1/r} \Big) = \Upsilon_r.
\end{align*}
Therefore, for $0 < r \leq \eps$ and $|\alpha | \leq \pi/2$:
\begin{displaymath}
\Re f_\theta (r \e^{i\alpha}) \geq (\cos \Upsilon_r) \,r \Big(\log \frac{1}{r}\Big)^\theta > 0,
\end{displaymath}
as announced in Lemma~\ref{definition f indice theta}. \qed
\bigskip

Let now $g_\theta$ be the conformal mapping from $\D$ onto $V_\eps$, which maps $\T = \partial \D$ 
onto $\partial V_\eps$, and with $g_\theta (1) = 0$ and $g'_\theta (1) = -\eps/4$. Explicitly, $g_\theta$ is 
the composition of the following maps: a) $\sigma \colon z \mapsto -z$ from $\D$ onto itself; b) 
$\gamma \colon z \mapsto \frac{z + i}{1 + iz} = \frac{ z + \bar{z} + i ( 1 - |z |^2)}{|1 + i z |^2}$ from $\D$ onto 
$P=\{ \Im z > 0\}$; c) $s \colon z \mapsto \sqrt z$ from $P$ onto $Q =\{\Re z > 0\,, \Im z >0\}$; d) 
$\gamma^{-1} \colon z \mapsto \frac{ z -i }{ 1 -iz}$ from $Q$ onto $V = \{ |z | < 1\,, \Re z >0\}$, and 
e) $h_\eps \colon z \mapsto \eps z$ from $V$ onto $V_\eps$.
\par\smallskip

Let:
\begin{equation}\label{definition phi indice theta}
\phi_\theta = \exp ( - f_\theta \circ g_\theta).
\end{equation}

By Lemma~\ref{definition f indice theta}, the analytic function $\phi_\theta$ maps $\D$ into 
$\D$. Moreover, one has $| \phi_\theta^\ast | < 1$ on $\partial \D \setminus \{1\}$.
\par\medskip

Now Theorem~\ref{theo Schatten} will follow from the following proposition.

\begin{proposition}\label{prop Schatten}
With the above notation, $C_{\phi_\theta} \colon H^2 \to H^2$ is compact for every $\theta >0$ and, 
moreover, is in the Schatten class $S_p$ if and only if $p > \frac{4}{\theta}\cdot$
\end{proposition}

Hence, given $p_0 > 0$, if we choose $\theta_0 = \frac{4}{p_0}$, we get that $C_{\phi_{\theta_0}}$ 
is in $S_p$ if and only if $p > p_0$. \qed
\par\medskip

\noindent{\bf Proof of Proposition~\ref{prop Schatten}} The compactness was proved in \cite{Shap-Tay}, 
Theorem~4.2; indeed, setting $h_\theta = f_\theta \circ g_\theta$, J. Shapiro and P. Taylor used the 
symbol $\phi_\theta = 1 - h_\theta$ and proved that $\lim_{t \to 0} | \phi_\theta '(\e^{it}) | = + \infty$:  
see equation (4.4) in \cite{Shap-Tay} (note that, in order to deduce the compactness from this equality, 
J. Shapiro and P. Taylor had to prove a theorem: Theorem~2.4 in \cite{Shap-Tay}, whose proof needs to use 
Gabriel's Theorem). Since our symbol is $\phi_\theta = \exp ( - h_\theta)$, and 
$\exp\big( h_\theta (\e^{it}) \big) \to 1$ 
as $t \to 0$, it follows that the derivatives have the same behaviour when $t \to 0$. However, we are 
going to recover this result by another method. For convenience, we shall write the boundary values 
$f_\theta^\ast, g_\theta^\ast, \ldots$ of the different analytic functions $f_\theta, g_\theta,\ldots$ in the 
same way as the analytic functions, without the exponent ${}^\ast$.\par
\smallskip

Note that, by Lemma~\ref{definition f indice theta}, $| \phi_\theta |$ is far from $1$ when 
$g_\theta (z)$ belongs to the half-circle $\{\eps \e^{i\alpha}\,;\ |\alpha | \leq \pi/2 \}$. Hence, we only have 
to study the case where $g_\theta ( \e^{i t}) = it$, $ -\eps \leq t \leq \eps$. Moreover, since 
$f_\theta (it) = it \big(\log \frac{1}{| t |} - i\, {\rm sgn}\,(t) \frac{\pi}{2}\big)^\theta$, one has
\begin{displaymath}
| f_\theta (it) | \approx | t |\Big(\log \frac{1}{| t |}\Big)^\theta,
\end{displaymath}
so that $| \phi_\theta \big( g_\theta^{-1}(it)\big) |$ is far from $1$ when $t$ is away from $0$. Therefore, 
it suffices to study what happens when $t$ is in a neighbourhood of $0$.\par
Note also that $g_\theta$ is bi-Lipschitz in a neighbourhood of $1$ (so $g_\theta (\e^{it}) \approx it$ 
when $|t |\leq \pi/2$), so we may forget it, and only consider the measure of the $t$'s for 
which $f_\theta (it)$ belongs to the suitable sets. Moreover, for convenience, we only write the proof for 
$t > 0$.\par
\smallskip

Since
\begin{align}
f_\theta (it) 
& = it \Big(\log \frac{1}{t} \Big)^\theta 
\bigg[\Big(1 - \frac{i\pi/2}{\log 1/t} \Big)^\theta \bigg]  \notag\\
& = it \Big(\log \frac{1}{t} \Big)^\theta \bigg[1 - \frac{i \pi \theta /2}{\log 1/t} 
+ o\,\Big(\frac{1}{\log 1/t}\Big)\bigg], \label{developpement f theta de it}
\end{align}
one has:
\begin{align}
& \Re f_\theta (it) \approx t \Big(\log \frac{1}{t}\Big)^{\theta - 1} \label{partie Re de f theta}\\
& \Im f_{\theta} (it) \approx t \Big(\log \frac{1}{t}\Big)^\theta. \label{partie Im de f theta}
\end{align}
\par
Now, we have $\exp\big( - f _\theta (it ) \big) \in W (e^{- i\alpha}, h)$ ($0\leq \alpha < 2\pi$) if and only if 
\begin{align}
& \Re f_\theta (it) \lesssim h  \label{numero un} \\
& | \Im f_{\theta} (it) - \alpha | \lesssim h. \label{numero deux}
\end{align}
But, when $e^{- i\alpha} \not \simeq 1$, this cannot happen for $h$ small enough (since $t$ goes to $0$ 
as $h$ goes to $0$). Essentially, we only have to consider the case $\e^{- i\alpha} = 1$, for which one has: 
when $t > 0$ goes to $0$, $|\Re f_\theta (it) | \lesssim |\Im f_{\theta} (it) |$, and hence 
$\exp\big( - f _\theta (it ) \big) \in W (1, h)$ if and only if $t (\log 1/t)^\theta \lesssim h$, \emph{i.e.} 
$t \lesssim h/ (\log 1/h)^\theta$. Actually, we cannot assume $\alpha = 0$, and we have to be more precise, 
and must do the following reasoning. When  \eqref{numero un} is satisfied, one has 
$t (\log 1/t)^{\theta - 1} \lesssim h$; then $t \lesssim h/(\log 1/h)^{\theta -1}$ and hence: 
\begin{displaymath}
t \Big( \log \frac{1}{t} \Big)^\theta \lesssim 
\frac{h}{(\log 1/h)^{\theta - 1}} \Big( \log \frac{1}{h}\Big)^\theta 
= h \log \frac{1}{h}\cdot
\end{displaymath}
It follows that the condition~\eqref{numero deux}  implies that:
\begin{equation}\label{encadrement alpha}
0 \leq \alpha \lesssim h \log \frac{1}{h}\, \cdot
\end{equation}
Therefore condition~\eqref{numero deux} implies that $t (\log 1/t)^\theta \lesssim h \log 1/h$, which gives: 
\begin{displaymath}
t \lesssim \frac{h \log 1/h}{(\log 1/h)^\theta} = \frac{h}{(\log 1/h)^{\theta - 1}} \raise 1,5 pt \hbox{,}
\end{displaymath}
\emph{i.e.} $t(\log 1/t)^{\theta - 1} \lesssim h$: condition~\eqref{numero un} is satisfied (up to a constant 
factor for $h$).\par
Since, by \eqref{encadrement alpha}, condition~\eqref{numero deux} is satisfied when 
$\alpha - h \lesssim t (\log 1/t)^\theta \lesssim (\alpha + h)$ and implies that 
$ - (\alpha + h) \lesssim t (\log 1/t)^\theta \lesssim (\alpha + h)$, it follows that the set of $t$'s 
such that \eqref{numero un} and \eqref{numero deux} are satisfied has, since 
$\log (\alpha + h) \approx \log h$, a measure $\approx h/(\log 1/h)^\theta$.
\par\smallskip

We then have proved that:
\begin{equation}\label{rho pour phi indice theta}
\rho_{\phi_\theta} (h) \approx \frac{h}{(\log 1/h )^\theta}\cdot
\end{equation}
Since $\rho_{\phi_\theta} (h) = o\, (h)$, MacCluer's criterion gives the compactness of $C_{\phi_\theta}$.
\medskip

Now, we shall examine when $C_{\phi_\theta}$ is in the Schatten class $S_p$, and for that we shall use 
Luecking's theorem (Theorem~\ref{theo Luecking}). We have to analyze the behaviour of $\phi_\theta$ with 
respect to the Luecking sets $R_{n,j}$; actually, for convenience, we shall work with 
$R'_{n,j} = R_{n, 2^n - 1 - j}$. \par
\smallskip

We have to consider the modulus and the argument constraints.
\par\smallskip

\emph{Modulus constraint.} The condition $\frac{h}{2} \leq \Re f_\theta (it) < h$ writes:
\begin{equation}\label{encadrement module f indice theta}
\frac{h}{2} \lesssim t \Big(\log \frac{1}{t}\Big)^{\theta - 1} \lesssim h,
\end{equation}
and reads as: 
\begin{equation}\label{encadrement module avec h}
\frac{h/2}{(\log 1/h)^{\theta - 1}} \lesssim t \lesssim \frac{h}{(\log 1/h)^{\theta - 1}}\,\raise 1,5 pt \hbox{,}
\end{equation}
or, when $h= h_n = 2^{-n}$:
\begin{equation}\label{contrainte module Luecking}
\frac{2^{-(n+1)}}{(n+1)^{\theta -1}} \lesssim t \lesssim \frac{2^{-n}}{n^{\theta - 1}}\cdot
\end{equation}
More precisely, one must have $t \in I_n = [a_n, b_n]$, with 
$a_n \approx \frac{2^{-(n+1)}}{(n+1)^{\theta -1}}$\raise 1,5 pt \hbox{,} 
$b_n \approx \frac{2^{-n}}{n^{\theta - 1}}$\raise 1,5 pt \hbox{,}
and $|I_n| \approx \frac{2^{-n}}{n^{\theta - 1}}\cdot$
\smallskip

\emph{Argument constraint.} We are looking for the set $J_n$ of the indices $j = 0,1, \ldots, 2^n - 1$ for 
which we have $\exp\big( - f_\theta (it) \big) \in R'_{n,j}$. We have must have both the modulus constraint 
$t \in I_n$ and:
\begin{equation}\label{contrainte argument theta}
j h_n \lesssim t \Big(\log \frac{1}{t} \Big)^\theta \lesssim (j+1) h_n, 
\end{equation}
which implies, for $j \geq 1$, since one has \eqref{encadrement module f indice theta}:
\begin{displaymath}
j \lesssim \log \frac{1}{t} \lesssim 2 (j + 1),
\end{displaymath}
and hence, by \eqref{encadrement module avec h} and \eqref{contrainte module Luecking}:
\begin{displaymath}
\frac{n}{2} \log 2 \lesssim j \lesssim n \log 2,
\end{displaymath}
\emph{i.e.} $j \approx n$. The constant coefficients are not relevant here, and hence this estimation means 
that $j$ can only take $O (n)$ values.\par
On the other hand, when the modulus constraint \eqref{encadrement module f indice theta} is satisfied,  one has 
$\log (1/t)  \approx \log (1/h_n)$, and the argument constraint~\eqref{contrainte argument theta} is equivalent, 
with $h=h_n$, to:
\begin{displaymath}
\frac{j h}{ (\log 1/h)^\theta} \lesssim t \lesssim \frac{(j+1)h}{(\log 1/h)^\theta} \cdot
\end{displaymath}
The length of the corresponding interval is $\approx h/ \big(\log (1/h)\big)^\theta $, which is equal, when 
$h = h_n = 2^{-n}$, to $2^{-n} / n^\theta$.\par
It follows hence that $m_{\phi_\theta} (R'_{n,j}) \approx 2^{-n} / n^\theta$ for exactly $O(n)$ values of $j$, 
and otherwise $m_{\phi_\theta} (R'_{n,j}) =0$; therefore:
\begin{displaymath}
\sum_{n=0}^\infty 2^{np/2} \sum_{j=0}^{2^n - 1} \big[m_{\phi_\theta} (R'_{n,j}) \big]^{p/2}
\approx \sum_{n=1}^\infty 2^{np/2} n \Big( \frac{2^{-n}}{n^\theta} \Big)^{p/2} 
= \sum_{n=1}^\infty \frac{1}{n^{\theta \frac{p}{2} - 1}} \raise 1,5pt \hbox{,}
\end{displaymath}
which is finite if and only if $\frac{\theta p}{2} - 1 >1$, \emph{i.e.} $p > 4/ \theta$. \qed
\par
\bigskip

\noindent{\bf Remark 1.} Proposition~\ref{prop Schatten} shows that, for these symbols $\phi_\theta$, the 
necessary condition~\eqref{limite condition necessaire Schatten} of 
Proposition~\ref{condition necessaire Schatten} is not sharp; in fact, we have proved in 
\eqref{rho pour phi indice theta} that $\rho_{\phi_\theta} (h) \approx h /(\log 1/h)^\theta$; hence 
\eqref{limite condition necessaire Schatten} gives $\theta > 2/p$ when $C_{\phi_\theta} \in S_p$, even 
though we must have $\theta > 4/p$.
\bigskip

\noindent{\bf Remark 2.} These composition operators answer negatively the question asked at the end 
of Section~\ref{section same modulus}, since, by \eqref{rho pour phi indice theta}, the measure 
$m_{\phi_\theta}$ is $\alpha$-Carleson for no $\alpha > 1$, though $C_{\phi_\theta}$ is 
in $S_p$ for every $p >  4/\theta$.\par
\bigskip

\noindent{\bf Remark 3.} Our proof of Theorem~\ref{meme module puissance} fails when 
$\alpha =2$ (\emph{i.e.} $\beta =1$). However, we are going to see that in this case, the composition 
operator $C_{\phi_1}$ of this Theorem~\ref{meme module puissance} is in $S_p$ for every $p>4$.\par
Indeed, for $\beta =1$, one has, in an explicit way:
\begin{displaymath}
\Big|\sin \frac{t}{2} \Big| =\frac{2}{\pi} - \frac{4}{\pi} \sum_{k=1}^\infty \frac{\cos kt}{4 k^2 - 1} 
=c_0 +\sum_{k=1}^\infty c_k \cos kt \,,
\end{displaymath}
with $c_k = O (k^{-2})$. Then: 

\begin{displaymath}
{\cal H} f (t) =\sum_{k=1}^\infty \frac{\sin kt}{4 k^2 - 1} = \frac{1}{4} u (t) + v (t),
\end{displaymath}
where:
\begin{displaymath}
u (t) = \sum_{k=1}^\infty \frac{\sin kt}{k^2}\,\raise 1,5pt \hbox{,}
\end{displaymath}
and $v$ is ${\cal C}^2$, since $\hat v (k) = O (k^{-4})$. We have:
\begin{displaymath}
\qquad \qquad u ' (t) = \sum_{k=1}^\infty \frac{\cos kt}{k} 
= - \log \Big(2 \sin \frac{t}{2}\Big)\,\raise 1,5pt \hbox{,} 
\qquad 0 < t < \pi,
\end{displaymath}
(note that $u \notin {\rm Lip}\,1$, but $u$ is in the Zygmund class). It follows, when $t > 0$ goes to zero, that:
\begin{displaymath}
{\cal H}f (t) \approx - t \log \Big( 2\sin \frac{t}{2} \Big) \approx t \log \frac{1}{t} \cdot
\end{displaymath}
Since the modulus of continuity of ${\cal H}f$ at $0$ is $\approx t \log (1/t)$ (see also \cite{Zyg}, Chapter III, 
Theorem~(13.30) for a general result), it follows that the argument condition 
on $\phi_1$ is $|t| \lesssim \frac{h}{\log (1/h)}$ and that 
$\rho_{\phi_1}(h) =O\big(\frac{h}{\log (1/h)}\big)$. Then $C_{\phi_1}$ is compact. It is actually in $S_p$ 
for every $p>4$. In fact, when $t > 0$ goes to zero, $\sin (t/2) \sim t/2$, and hence, with the notation of the 
Subsection~\ref{general construction}, one has $| \Phi^\ast (t) | \sim t$, whereas 
$\arg \big( \Phi^\ast (t) \big) \approx t/ \log (1/ t)$. We have then the same conditions as in the proof of 
Proposition~\ref{prop Schatten}, when $\theta =1$. Hence $C_{\phi_1} \in S_p$ for every $p>4$. \qed
\bigskip

\noindent{\bf Remark 4.} We shall characterize in \cite{Orlicz-Hardy} the Orlicz functions $\Psi$ for which 
these composition operators $C_{\phi_\theta}$ are compact.
\bigskip

In Theorem~\ref{theo Schatten}, we get, for any $p_0 > 0$, a composition operator which is in $S_p$ if and 
only if $p > p_0$. We can modify slightly this operator so as to belong in $S_p$ if and only if $p \geq p_0$. 

\begin{theoreme}\label{theo Schatten bis}
For every $p_0 >0$, there exists an analytic self-map $\phi \colon \D \to\D$ such that the composition 
operator $C_\phi \colon H^2 \to H^2$ is in the Schatten class $S_{p_0}$, but not in $S_p$, for $p < p_0$.
\end{theoreme}

\noindent{\bf Proof.} We shall use the same method as in Theorem~\ref{theo Schatten}, but by replacing 
the function $f_\theta$ by this modified function:
\begin{equation}
\tilde f_\theta (z) = z (- \log z )^\theta \big[ \log ( - \log z) \big]^q,
\end{equation}
where $q > \theta /2$ ($q =\theta$, for example).\par
Call $\tilde \phi_\theta$ the corresponding self-map of $\D$.\par
\smallskip

We shall not give all the details since they are mostly the same as in the proof of 
Theorem~\ref{theo Schatten}.\par
We have first to check that:
\begin{lemme}\label{petit lemme}
For $\eps >0$ small enough, one has  $\Re \tilde f_\theta (z) > 0$ for every $z \in V_\eps$.
Moreover $|\tilde \phi_\theta^\ast | <  1$ a.e. on $\partial \D$.
\end{lemme}

\noindent{\bf Proof.} Write $z = r \e^{i\alpha}$ with $0 < r < \eps$ and $|\alpha | \leq \pi/2$. One has:
\begin{displaymath}
\tilde f_\theta (z) = r \e^{i\alpha} \Big( \log \frac{1}{r} - i \alpha\Big)^\theta 
\bigg( \log \Big( \log \frac{1}{r} - i \alpha \Big)\bigg)^q.
\end{displaymath}
But:
\begin{displaymath}
\Big( \log \frac{1}{r} - i \alpha\Big)^\theta = \Big( \log \frac{1}{r}\Big)^\theta 
\bigg( 1 - \frac{i \alpha \theta}{\log 1/r} + o\, \Big( \frac{1}{( \log 1/r) }\Big) \bigg);
\end{displaymath}
and, on the other hand:
\begin{align}
\log \Big( \log \frac{1}{r} - i \alpha \Big) 
& = \log \sqrt{\Big( \log \frac{1}{r} \Big)^2 + \alpha^2 } 
+ i \arg \Big(  \log \frac{1}{r} - i \alpha \Big)  \label{log log}\\
& = \Big(\log \log \frac{1}{r} \Big) \big( 1 + o\,(1/ \log (1/r)) \big) 
- i \arctan \frac{\alpha}{\log 1/r} \notag \\
& = \Big( \log \log \frac{1}{r} \Big) \big( 1 + o\,(1/ \log (1/r) ) \big)
- i \frac{\alpha}{\log 1/r} \big( 1 + o\,(1) \big)\,\raise 1,5 pt \hbox{,} \notag
\end{align}
and
\begin{align}\label{puissance q}
\bigg[& \log \Big( \log \frac{1}{r}  - i \alpha \Big) \bigg]^q \\
& = \Big(\log\log \frac{1}{r} \Big)^q \bigg[ \Big[ 1 + o\,\Big(\frac{1}{\log 1/r}\Big)\Big] 
- i \frac{1}{\log 1/r \log\log 1/r} \big( \alpha q + o\,(1)\big) \bigg]  \notag \\
& =  \Big(\log\log \frac{1}{r} \Big)^q \bigg[ 1 + o\,\Big(\frac{1}{\log 1/r}\Big) \bigg]. \notag
\end{align}
Hence, for $\eps >0$ small enough: 
\begin{align*}
\Re \tilde f_\theta (z) 
& = r (\log 1/r )^\theta ( \log \log 1/r)^q 
\bigg[ \cos \alpha +\frac{\theta \alpha \sin \alpha}{\log 1/r} + o\,\Big( \frac{1}{\log 1/r} \Big) \bigg] \\
& \geq r (\log 1/r )^\theta ( \log \log 1/r)^q \bigg[ \cos \alpha 
+ \frac{\theta (\alpha \sin \alpha - 1/4 ) }{\log 1/r} \bigg].
\end{align*}
That gives the result since, on the one hand, $|\theta (\alpha \sin \alpha - 1/4 ) / \log 1/r |$ is $\leq \sqrt 2/4$ 
for $\eps > 0$ small enough and $\cos \alpha \geq \sqrt 2/2$ when $|\alpha | \leq \pi/4$, and, on the 
other hand, when $|\alpha | \geq \pi/4$, one has 
\begin{displaymath}
\cos \alpha + \theta (\alpha \sin \alpha - 1/4 ) \geq \theta \big( \pi \sqrt 2/8 - 1/4) > 0.
\end{displaymath}
 That ends the proof of the lemma. \qed
\medskip

Now, by \eqref{log log}, for $t > 0$ going to zero:
\begin{displaymath}
\log \Big( \log \frac{1}{t} - i \frac{\pi}{2} \Big) 
= \bigg[ \log \log \frac{1}{t} + O \Big( \frac{1}{(\log t)^2}\Big) \bigg] + i\, O \Big(\frac{1}{\log 1/t}\Big)\cdot
\end{displaymath}

Therefore:
\begin{align*}
& \Re \bigg[ \log \Big( \log \frac{1}{t} - i \frac{\pi}{2} \Big) \bigg]^q 
\approx \Big( \log \log \frac{1}{t} \Big)^q \\
& \Im  \bigg[ \log \Big( \log \frac{1}{t} - i \frac{\pi}{2} \Big) \bigg]^q
\approx \frac{(\log \log 1/t)^{q - 1}}{\log 1/t}\, \raise 1,5pt \hbox{,}
\end{align*}
and, using \eqref{developpement f theta de it} (or \eqref{puissance q}, with 
$\alpha = \pi/2$), we get:
\begin{align*}
& \Re \tilde f_\theta (it) \approx t \Big( \log \frac{1}{t} \Big)^\theta 
\frac{(\log \log 1/t)^q}{\log 1/t} 
= t \Big( \log \frac{1}{t} \Big)^{\theta - 1} \Big(\log \log \frac{1}{t} \Big)^q  \\
& \Im \tilde f_\theta (it) \approx t \Big( \log \frac{1}{t} \Big)^\theta 
\Big( \log \log \frac{1}{t} \Big)^q.
\end{align*}
Hence the modulus constraint gives, with $h = h_n = 2^{-n}$:
\begin{equation}\label{blabla}
t \approx \frac{h}{\big( \log 1/h \big)^{\theta - 1} \big( \log \log 1/h\big)^{q - 1}}\,\raise 1,5 pt \hbox{,}
\end{equation}
and the argument constraint:
\begin{displaymath}
jh \lesssim t \Big( \log \frac{1}{t}\Big)^\theta \Big( \log \log \frac{1}{t} \Big)^q \lesssim (j+1)h.
\end{displaymath}
One gets:
\begin{displaymath}
j \approx \log\frac{1}{t} \approx \log \frac{1}{h} \approx n,
\end{displaymath}
and:
\begin{displaymath}
m_{\tilde \phi_\theta} (R'_{n,j}) \approx \frac{h}{\big( \log 1/h)^\theta \big( \log \log 1/h\big)^q} 
\approx \frac{2^{-n}}{n^\theta (\log n)^q}\cdot
\end{displaymath}
It follows that Luecking's criterion becomes:
\begin{displaymath}
\sum_{n=1}^\infty 2^{np/2} n \frac{2^{-np/2}}{n^{\theta p/2} (\log n)^{q p/2}} 
= \sum_{n=1}^\infty \frac{1}{n^{\theta p/2 - 1} (\log n)^{q p/2}}\cdot
\end{displaymath}
\par

This series converges if and only if $p > 4/\theta$ or else $p= 4/\theta$, since then 
$q p/2 = q (4/\theta) /2 = q/ (\theta /2) > 1$. Hence $C_{\tilde \phi_\theta} \in S_p$ if and only if 
$p \geq 4/\theta$, and that proves Theorem~\ref{theo Schatten bis}. \qed

\bigskip

In \cite{Carroll-Cowen}, T. Carroll and C. Cowen showed that there exist compact composition operators on 
$H^2$ which are in no Schatten class $S_p$ for $p < \infty$ (see also \cite{Zhu} and \cite{Jones}). We shall 
give an explicit example of such an operator.

\begin{theoreme}\label{no Schatten}
There exist compact composition operators $C_\phi \colon H^2 \to H^2$ which are in no Schatten class 
$S_p$ with $p < \infty$.
\end{theoreme}

\noindent{\bf Proof.} The proof follows the lines of those of Theorem~\ref{theo Schatten} and 
Theorem~\ref{theo Schatten bis}, but with, instead of $f_\theta$ or $\tilde f_\theta$:
\begin{equation}
f ( z) = z \log (- \log z).
\end{equation}
For $\eps >0$ small enough, we have $\Re f (z) > 0$, so that the corresponding function $\phi$ sends 
$\D$ into itself; moreover, one has $|\phi^\ast | < 1 $ \emph{a.e.} on $\partial \D$. Indeed, if 
$z = r \e^{i\alpha} \in V_\eps$, then $- \log z = R\,\e^{i\beta}$, with
\begin{align*}
R & = \sqrt{(\log 1/r)^2 +\alpha^2} = \log (1/r) + o\,(\log 1/r) \\
\beta & = - \arctan \frac{\alpha}{\log 1/r} = - \frac{\alpha}{\log 1/r} + o\, \big(1/ \log (1/r) \big).
\end{align*}
Then $\log (- \log z) = \log R + i\beta = (\log R) \big( 1 + \frac{i \beta}{\log R}\big)$ and 
\begin{displaymath}
\Re f (z) = r (\log R) \Big( \cos \alpha - \frac{\beta \sin \alpha}{\log R} \Big).
\end{displaymath}
But
\begin{displaymath}
\cos \alpha - \frac{\beta \sin \alpha}{\log R} = \cos \alpha + 
\frac{\alpha \sin \alpha}{(\log 1/r ) (\log \log 1/r) } + o\, \big(1/ (\log 1/r) (\log \log 1/r) \big)
\end{displaymath}
and we see, as in the proof of Lemma~\ref{petit lemme} that this quantity is, for $\eps >0$ small enough, 
greater than a positive constant, uniformly for $|\alpha | \leq \pi/2$.
\par
\smallskip

Now, one has, as above, for $t $ going to zero:
\begin{align*}
& \Re f (it ) \approx \frac{|t|}{\log 1/|t|} \\
& \Im f (it) \approx t \log \log \frac{1}{|t|}\cdot
\end{align*}
It follows that $C_\phi$ is compact since $\rho_\phi (h) \approx h/\log\log 1/h$: indeed, we have to control 
the two conditions:
\begin{align}
0 & < \frac{|t|}{\log 1/|t|} \lesssim h \label{un}\\
\alpha - h & \lesssim t \log \log \frac{1}{|t|}  \lesssim \alpha + h 
\,\raise 1,5pt \hbox{,} \label{deux}
\end{align}
for $0 \leq \alpha < 2\pi$. When $\alpha \not\simeq 0$, condition~\eqref{deux} cannot happen 
for $h$ small enough, because of condition~\eqref{un}; more precisely, \eqref{un} and \eqref{deux} 
imply that
\begin{equation}\label{trois}
0 \leq \alpha \lesssim h + h \log \frac{1}{h} \log \log \frac{1}{h} 
\lesssim  h \log \frac{1}{h} \log \log \frac{1}{h} \raise 1,5 pt \hbox{,}
\end{equation}
and then \eqref{deux} implies:
\begin{displaymath}
- h \lesssim t \log \log \frac{1}{|t|} \lesssim h \log \frac{1}{h} \log \log \frac{1}{h} \raise 1,5 pt \hbox{;}
\end{displaymath}
\emph{ a fortiori}:
\begin{displaymath}
|t|  \log \log \frac{1}{|t|} \lesssim h \log \frac{1}{h} \log \log \frac{1}{h} \raise 1,5 pt \hbox{,}
\end{displaymath}
and then:
\begin{displaymath}
|t| \lesssim h \log \frac{1}{h} \cdot
\end{displaymath}
Therefore condition~\eqref{deux} implies condition~\eqref{un}. Since, \eqref{trois} implies that 
\begin{displaymath}
\log (\alpha + h) \approx \log h , 
\end{displaymath}
we get that the measure of the $t$'s satisfying \eqref{deux} is $\approx h\log \log 1/h$: \eqref{deux}  is 
implied by $-h \leq t \log \log 1/|t| \leq \alpha +h$ and implies that 
$- (\alpha +h) \leq t \log \log 1/|t| \leq (\alpha +h)$.\par

Hence, we have got that $\rho_\phi (h) \approx h / \log \log 1/h$, showing that $C_\phi$ is compact.
\par
By Proposition~\ref{condition necessaire Schatten}, $C_\phi$ is in no Schatten class $S_p$, $p > 0$.\par
In order to make more transparent the behaviour of the measure $m_\phi$, we shall give a direct 
proof, using Luecking's criterion.
\smallskip

For testing Luecking's criterion, we must have, assuming $t > 0$:
\begin{displaymath}
t \approx h\log \frac{1}{h}
\end{displaymath}
and:
\begin{displaymath}
jh \lesssim t \log \log \frac{1}{t} \lesssim (j + 1) h,
\end{displaymath}
so $j \approx \log 1/h \log \log 1/h \approx n \log n$, and:
\begin{displaymath}
m_\phi (R'_{n,j}) \approx \frac{h}{\log \log 1/h} \approx \frac{2^{-n}}{\log n}\cdot
\end{displaymath}
Therefore Luecking's condition becomes:
\begin{displaymath}
\sum_{n=1}^\infty 2^{np/2} \bigg[n\log n \bigg(\frac{2^{-np/2}}{(\log n)^{p/2}} \bigg) \bigg]= 
\sum_{n=1}^\infty \frac{n}{(\log n)^{p/2 - 1}} = +\infty.
\end{displaymath}
It follows that $C_\phi \notin S_p$. \qed
\bigskip

As a corollary, we get:

\begin{theoreme}\label{no Schatten meme module}
There exist two symbols $\phi$ and $\psi$ of same modulus $|\phi^\ast| = |\psi^\ast|$ on $\T$ such that 
$C_\phi$ is compact on $H^2$, but in no Schatten class $S_p$ for $p < \infty$, whereas $C_\psi$ is in 
$S_p$ for every $p > 2$.
\end{theoreme}

\noindent{\bf Proof.} Set 
\begin{displaymath}
\psi = \phi M,
\end{displaymath}
where $\phi$ is the function used in Theorem~\ref{no Schatten}, and $M$ is the singular inner function given 
by \eqref{M}. Since $\Re f (it ) \approx t/ \log (1/t)$, Lemma~\ref{deuxieme lemme} gives, if one sets 
$\tilde f (t) = t/ \log (1/t)$:
\begin{displaymath}
\rho_\psi (h) \approx h \tilde f^{-1} (h) \approx h^2 \log (1/h) .
\end{displaymath}
It follows that $m_\psi$ is an $\alpha$-Carleson measure for any $\alpha <  2$. Hence, given any $p >2$, 
and choosing an $\alpha < 2$ such that $p > 2/(\alpha - 1)$, it follows from Corollary~\ref{coro Luecking} 
that $C_\psi \in S_p$. \qed
\bigskip

A funny question about Schatten classes is: for which Orlicz functions $\Psi$ is the above composition operator 
$C_\phi$ in the Orlicz-Schatten class $S_\Psi$?

%%%%%%%%%%%%%%%%
\section{Questions}

We shall end this paper with some comments and questions.\par
\smallskip

{\bf 1.} G. Pisier recently suggested that one should study more carefully the approximation numbers of 
composition operators. Since an operator $T$ on a Hilbert space is in the Schatten class $S_p$ if and 
only if $\sum_{n=1}^\infty a_n^p < +\infty$, where $a_n$  is the $n^{th}$ approximation number of $T$, 
our present work may be seen as a, very partial, contribution to that study. \par
Note that Luecking's proof of its trace-class theorem does not make explicit mention of singular numbers, 
but relies instead on an interpolation argument.
\par\smallskip

{\bf 2.} It is clear, by our results, that the membership in $S_p$ for the composition operator $C_\phi$ cannot 
be characterized in terms of the growth of the Carleson function $\rho_\phi$. But we have given a sufficient 
condition in Corollary~\ref{coro Luecking}, and a necessary one in 
Proposition~\ref{condition necessaire Schatten}. Are these two conditions sharp? Can they be improved?\par
\smallskip

{\bf 3.} Do there exist two symbols $\phi_1$ and $\phi_2$ having the same modulus on $\partial \D$ such 
that $C_{\phi_1}$ is not compact on $H^2$, but $C_{\phi_2}$ is in $S_p$ for every $p > 2$?\par
\smallskip

{\bf 4.} (K. Kellay) If $C_\phi \colon H^2 \to H^2$ is compact (or even in $S_p$, with $p >2$) and 
$\psi \colon \D \to \D$ is analytic, is $C_{\phi\psi} \colon H^2 \to H^2$ compact?\par
\smallskip

%%%%%%%%%%%%%%%%%%%%%%%%

\bigskip

\vbox{\noindent{\it 
{\rm Pascal Lef\`evre}, Universit\'e d'Artois,\\
Laboratoire de Math\'ematiques de Lens EA~2462, \\
F\'ed\'eration CNRS Nord-Pas-de-Calais FR~2956, \\
Facult\'e des Sciences Jean Perrin,\\
Rue Jean Souvraz, S.P.\kern 1mm 18,\\ 
62\kern 1mm 307 LENS Cedex,
FRANCE \\ 
pascal.lefevre@euler.univ-artois.fr 
\smallskip

\noindent
{\rm Daniel Li}, Universit\'e d'Artois,\\
Laboratoire de Math\'ematiques de Lens EA~2462, \\
F\'ed\'eration CNRS Nord-Pas-de-Calais FR~2956, \\
Facult\'e des Sciences Jean Perrin,\\
Rue Jean Souvraz, S.P.\kern 1mm 18,\\ 
62\kern 1mm 307 LENS Cedex,
FRANCE \\ 
daniel.li@euler.univ-artois.fr
\smallskip

\noindent
{\rm Herv\'e Queff\'elec},
Universit\'e des Sciences et Technologies de Lille, \\
Labo\-ratoire Paul Painlev\'e U.M.R. CNRS 8524, \\
U.F.R. de Math\'ematiques,\\
59\kern 1mm 655 VILLENEUVE D'ASCQ Cedex, 
FRANCE \\ 
queff@math.univ-lille1.fr
\smallskip

\noindent
{\rm Luis Rodr{\'\i}guez-Piazza}, Universidad de Sevilla, \\
Facultad de Matem\'aticas, Departamento de An\'alisis Matem\'atico,\\ 
Apartado de Correos 1160,\\
41\kern 1mm 080 SEVILLA, SPAIN \\ 
piazza@us.es\par}
}

\end{document}